\documentclass{article}
\usepackage{amsmath,amsthm,amsfonts,eucal}

\numberwithin{equation}{section}

\def\ca{{\mathcal A}}
\def\cb{{\mathcal B}}

\def\cd{{\mathcal D}}

\def\cf{{\mathcal F}}

\def\ch{{\mathcal H}}

\def\ck{{\mathcal K}}
\def\cl{{\mathcal L}}
\def\cam{{\mathcal M}}

\def\car{{\mathcal R}}

\def\bc{{\mathbb C}}

\def\bn{{\mathbb N}}
\def\br{{\mathbb R}}

\def\a{\alpha}
\def\b{\beta}
\def\g{\gamma}        
        \def\D{\Delta}
\def\eps{\varepsilon}

\def\th{\vartheta}

\def\l{\lambda}       \def\La{\Lambda}
\def\m{\mu}

\def\s{\sigma}
    
\def\t{\tau}

\def\f{\varphi}

        \def\O{\Omega}

\def\imply{\Rightarrow}

\def\ov{\overline}
\def\e#1{{\rm e}^{#1}}

\def\limt{\liminf_{t\to0^+}}
\def\limx{\liminf_{x\to\infty}}

\def\supp{{\rm supp}\ }

\def\ad#1{d_\infty(#1)}

\def\itm#1{\item[($#1$)]}

\def\mt{\ov\cam}

\def\pesci{\bowtie}
\def\re{{\rm Re\;}}

\def\Ainv{\cam}
\def\Aloc{\ca}

\def\ozero#1{{\rm ord}_0(#1)}
\def\oinf#1{{\rm ord}_\infty(#1)}

\newtheorem{Thm}{Theorem}[section]
\newtheorem{Cor}[Thm]{Corollary}
\newtheorem{Prop}[Thm]{Proposition}
\newtheorem{Lemma}[Thm]{Lemma}
\theoremstyle{definition}
\newtheorem{Dfn}[Thm]{Definition}
\newtheorem{exmp}[Thm]{Example}

\theoremstyle{remark}
\newtheorem{rem}[Thm]{Remark} 
 
\title{\huge  Singular traces, dimensions,\\
 and Novikov-Shubin invariants}
\author{Daniele Guido, Tommaso Isola\\
Dipartimento di Matematica,\\ Universit\`a di Roma ``Tor
Vergata'',\\ I--00133 Roma, Italy.}

\date{}
\begin{document}
\maketitle
\markboth{Asymptotic dimension package}
{Singular traces and Novikov-Shubin invariants}
\renewcommand{\sectionmark}[1]{}

%

 \setcounter{section}{-1}
 \section{Introduction.}\label{sec:intro}

In Alain Connes noncommutative geometry, integration of a function
$a$ on a noncommutative manifold is described via the formula $\int
a=\tau(a|D|^{-d})$, where $D$ is an unbounded selfadjoint operator
playing the role of the Dirac operator, $d$ is the dimension of the
manifold, and $\tau$ is a singular trace introduced by Dixmier and
corresponding to logarithmic divergences.

As a consequence the question of the existence of a non-trivial
integral can be described in terms of the singular traceability of the
compact operator $|D|^{-d}$, namely of the existence of a finite
non-trivial singular trace on the ideal generated by $|D|^{-d}$.

In this paper we show that, under suitable regularity conditions on
the eigenvalue sequence of $|D|$, the dimension $d$ can be uniquely
determined by imposing that $|D|^{-d}$ is singularly traceable, thus
providing a geometric measure theoretic definition for $d$. We remark
that with this choice of $d$ the operator $|D|^{-d}$ always produces a
singular trace which gives rise to a non-trivial integration
procedure, even when the logarithmic trace fails.

Singular traceability of compact operators has been described in
\cite{AGPS} via their eccentricity. We propose here another sufficient
condition in terms of the polynomial order of (the eigenvalue sequence
of) the given compact operator. More precisely, we prove that if $T$
has polynomial order equal to 1 it is singularly traceable. Since the
polynomial order of $|D|^{-d}$ is equal to $d$ times the polynomial
order of $|D|^{-1}$, this shows that the geometric measure theoretic
dimension $d$ coincides with the inverse of the polynomial order of
$|D|^{-1}$.

In the second part of the paper we discuss large scale counterparts of
the geometric measure theoretic dimension for the case of covering
manifolds.

Le us recall that in the case of von Neumann algebras with a
continuous trace singular traces can be constructed on the ideal
generated by a trace-compact operator $a$ when the eigenvalue function
$\mu_a(t)$ has a suitable asymptotics when $t\to\infty$, as in the
discrete case.

However, a new family of singular traces arises in this case,
detecting the asymptotic behaviour of $\mu_a(t)$ when $t\to0$
\cite{GI2}.

We prove here that if the polynomial order at $0$, $\ozero{a}$, is 1, then
$a$ is singularly traceable.

In the case of a non-compact covering $M$ of a compact Riemannian
manifold, Atiyah considered the von Neumann algebra with continuous
trace consisting of the bounded operators on $L^2(M)$ commuting with
the action of the covering group \cite{Atiyah}. In this case both
behaviours at $0$ and at $\infty$ of trace-compact operators may give
rise to singular traces.

As far as the behaviour at infinity is concerned, we note that the
dimension $n$ of $M$ is twice the inverse of the order at infinity of
$\Delta_k^{-1}$, therefore, by the result mentioned above,
$\Delta_k^{-n/2}$ is singularly traceable at $\infty$, and the
dimension of $M$ has a geometric measure theoretic interpretation.

We prove here that the large-scale counterpart of this dimension,
namely the number $d=2\ozero{\Delta_k^{-1}}^{-1}$, coincides with the
$k$-th Novikov-Shubin number.

As a consequence, Novikov-Shubin numbers can be considered as
(asymptotic) dimensions in the sense of geometric measure theory,
since the corresponding power of $\Delta_k^{-1/2}$ produces a singular
trace. Let us remark that an asymptotic-dimensional interpretation of
the $0$th Novikov-Shubin number is given by the fact, proved by
Varopoulos \cite{Varopoulos1}, that $\alpha_0$ coincides with the growth of
the covering group.

Finally we show that the Novikov-Shubin number $\underline{\alpha}_k$
coincides with (the supremum of) the dimension at $\infty$ of the
semigroup $e^{-t\Delta_k}$ introduced by Varopoulos, Saloff-Coste and
Coulhon \cite{VSC}.

Generalizations of the mentioned results on covering manifolds to the
case of open manifolds with bounded geometry are studied in
\cite{GI2}.


 \section{Polynomial order of an operator}\label{sec:singular}

 In this section we give a quick review to some notions pertaining to
singular traces \cite{GI1}, and then give some sufficient conditions
for the construction of singular traces in terms of the local or
asymptotic polynomial order of an operator.

Let $\cam\subset \cb(\ch)$ be a  semifinite von Neumann algebra equipped with a
normal semifinite faithful trace $tr$. We refer the reader to \cite{Ta} for the
general theory of von Neumann algebras.  Let $\ov \cam$ be the collection of the
closed, densely defined operators $x$ on $\ch$ affiliated with $\cam$ such that
$tr(e_{|x|}(t,\infty))<\infty$ for some $t>0$.
 
$\ov \cam$, equipped with strong sense operations \cite{Se} and with the topology
of convergence in measure (\cite{Stn}, \cite{Ne}), becomes a topological
$*$-algebra, called the algebra of $tr$-measurable operators.

For any $x=\int_0^\infty t\ de_x(t)\in\ov \cam_+$, $E\in\br\to\nu_x(E):=tr(e_x(E))$ 
is a Borel measure on $\br$, and $tr(x):=\int_0^\infty t\ d\nu_x(t)$ is a
faithful extension of $tr$ to $\ov \cam_+$.

 \begin{Dfn} \label{1.1} {\rm \cite{FK}} Let $a\in\ov \cam$, and
define, for all $t>0$,
 \itm{i} $\l_a(t):=tr(e_{|a|}(t,\infty))$, the
distribution function of $a$ w.r.t. $tr$, 
 \itm{ii} $\m_a(t):=\inf\{s\geq0: \l_a(s)\leq t\}$, the
non-increasing rearrangement of $a$ w.r.t. $tr$, which is a
non-increasing and right-continuous function. Moreover
$\lim_{t\downarrow0}\mu_a(t)=\|a\|\in[0,\infty]$. 
 \end{Dfn}

 \begin{rem}\label{1.2} 
	 $(i)$ If $\cam:=L^\infty(X,m)$ and $tr(f):=\int fdm$, then $\mt$ is 
	 the $^*$-algebra of functions that are bounded except on a set of 
	 finite $m$-measure, and, for any $f\in\mt$, $\mu_f\equiv f^*$  
	 is the classical non-increasing rearrangement of $f$ \cite{BS}.  \\
	 $(ii)$ If $\cam=\cb(\ch)$ and $tr$ is the usual trace, then 
	 $\mu_a=\sum_{n=0}^\infty s_n\chi_{[n,n+1)}$, where $\{s_n\}$ is 
	 the sequence of singular values of the operator $a$, arranged in 
	 non-increasing order and counted with multiplicity \cite{Si}.
 \end{rem}

 \begin{Dfn}\label{1.3} {\rm \cite{GI1}} 
	 Let $a\in\mt$. Then
	 \itm{i} $a$ is called $0$-eccentric if, setting
	 $$
	 S^0_a(t):=
	 \begin{cases}
		 \int_0^t \m_a(s)ds & \m_a\in L^1(0,1)\\
		 \vbox{\vskip 0.2cm}\\
		 \int_t^1 \m_a(s)ds & \m_a\not\in L^1(0,1),\\
	 \end{cases}
	 $$
	 $1$ is a limit point of $\{ \frac{S^0_a(2t)}{ S^0_a(t)} \}$, when 
	 $t\to0$.  \\
	 \itm{ii} $a$ is called $\infty$-eccentric if, setting
	 $$
	 S^\infty_a(t):=
	 \begin{cases}
		 \int_1^t \m_a(s)ds & \m_a\not\in L^1(1,\infty) \\
		 \vbox{\vskip 0.2cm} \\
		 \int_t^\infty \m_a(s)ds & \m_a\in L^1(1,\infty), \\
	 \end{cases}
	 $$
	 $1$ is a limit point of $\{ \frac{S^\infty_a(2t)}{ S^\infty_a(t)} 
	 \}$, when $t\to\infty$.  \\
 \end{Dfn}

 In $(i)$ we could replace $(0,1)$ with $(0,c)$ for any $c<\infty$, 
 and in $(ii)$ we could replace $(1,\infty)$ with $(c,\infty)$ for any 
 $c<\infty$, without affecting the definitions.

 \begin{Thm}\label{1.4} {\rm (\cite{GI1}, Theorem 6.4) } 
	 Let $(\cam,tr)$ be a semifinite von~Neumann algebra with a normal 
	 semifinite faithful trace, $a\in \ov \cam_+$ an eccentric operator.  
	 Then there exists a singular trace $\t$, namely a trace vanishing 
	 on projections which are finite w.r.t. $tr$, whose domain is the 
	 measurable bimodule generated by $a$ and such that $\t(a)=1$.
 \end{Thm}
 
 Now we give sufficient conditions to ensure eccentricity at $0$, or at 
 $\infty$.  It is based on the notions of order of infinite at $0$ and 
 of order of infinitesimal at $\infty$.

 \begin{Dfn} 
	 For $a\in\mt$ we define
	 \itm{i} order of infinitesimal of $a$ at $\infty$ 
	 $$
	 \oinf{a} := \liminf_{t\to\infty} \frac{\log \m_{a}(t)}{\log 
	 (1/t)}\; ,
	 $$
	 \itm{ii} order of infinite of $a$ at $0$
	 $$
	 \ozero{a} := \liminf_{t\to0} \frac{\log \m_{a}(t)}{\log (1/t)}\; .
	 $$
 \end{Dfn}
 
 \begin{rem}\label{rem1.6} 
	 If $a\in\mt_{+}$, then for any $\a>0$, $\oinf{a^{\a}}= \a\ 
	 \oinf{a}$ and $\ozero{a^{\a}}= \a\ \ozero{a}$.
 \end{rem}

 \begin{Thm}\label{1.5} 
	 Let $a\in\mt$ be s.t. $\ozero{a}=1$.  Then $a$ is $0$-eccentric.
 \end{Thm}
 \begin{proof}
	 It is a consequence of the following Propositions.
 \end{proof}

 \begin{Prop}\label{1.6} 
	 Suppose $\m_a\not\in L^1(0,1)$, and $\limt 
	 \frac{\log\m_a(t)}{\log(1/t)}=1$.  \\
	 Then $\limsup_{t\to0^+}\frac{S^0_a(2t)}{ S^0_a(t)}=1$.
 \end{Prop}
 \begin{proof}
	 Let us set $x:=\log(1/t)$, and $g(x):=t\m_a(t)$.  Then 
	 $\frac{\log\m_a(t)}{\log(1/t)}=\frac{\log g(x)}{ x} +1$, so that 
	 the hypothesis becomes $\limx \frac{\log g(x)}{ x}=0$.  Besides, 
	 as $S(t)\equiv S^0_a(t) = \int_t^1 \m_a(s)ds$, $S$ is 
	 nonincreasing and convex, so that $S(t) \geq S(2t) \geq S(t) - 
	 t\m_a(t)$, that is $1 \geq \frac{S(2t)}{ S(t)} \geq 
	 1-\frac{t\m_a(t) }{ S(t)}$, and the thesis is implied by $\limt 
	 \frac{t\m_a(t) }{ S(t)}=0$, but, as $\frac{t\m_a(t) }{ S(t)} = 
	 \frac{g(x) }{ \int_0^x g(s)ds}$, setting $G(x):= \int_0^x 
	 g(s)ds$, the thesis follows from $\limx 
	 \frac{g(x)}{ G(x)}=0$.  \\
	 So we have to prove that $\limx \frac{\log g(x)}{ x}=0$ implies 
	 $\limx \frac{g(x)}{ G(x)}=0$.  Suppose on the contrary that there 
	 are $x_0,\eps>0$ s.t. $\frac{g(x)}{ G(x)}\geq\eps$, for $x\geq 
	 x_0$.  Then, for $x\geq x_0$, we have $\log\frac{G(x)}{ 
	 G(x_0)}=\int_{x_0}^x \frac{G'(s)}{ G(s)}ds\geq \eps(x-x_0)$, 
	 which implies $G(x)\geq G(x_0)\e{\eps(x-x_0)}$, and $g(x)\geq 
	 \eps G(x) \geq \eps G(x_0) \e{\eps(x-x_0)}$, so that $\frac{\log 
	 g(x) }{ x} \geq \frac{k }{ x}+\eps$ and $\limx \frac{\log g(x) }{ 
	 x}\geq \eps>0$, which is absurd.
 \end{proof}

 \begin{Prop}\label{1.7} 
	 Suppose $\m_a\in L^1(0,1)$ and $\limt 
	 \frac{\log\m_a(t)}{\log(1/t)}=1$.  \\
	 Then $\liminf_{t\to0^+}\frac{S^0_a(2t)}{ S^0_a(t)}=1$.
 \end{Prop}
 \begin{proof}
	 Recall that $S(t)\equiv S^0_a(t)=\int_0^t \m_a(s)ds$, is 
	 positive, nondecreasing, concave, and $S(0)=0$, $S(1)=1$, as we 
	 can multiply $\m_a$ by a suitable positive constant without 
	 altering our statement.  \\
	 Then $\log S$ is nondecreasing and concave, and 
	 $\lim_{t\to0^+}\log S=-\infty$, so that $g:=\frac{S'}{ S}$ is 
	 positive, nonincreasing, $\lim_{t\to0^+}g(t)=\infty$, $\log S(t) 
	 = -\int_t^1 g(s)ds$, and $\int_0^1 g(s)ds = -\lim_{t\to0^+}\log S 
	 = \infty$.  Finally $\m_a=S'=gS$, so that 
	 $\frac{\log\m_a(t)}{\log(1/t)} = \frac{\log(tg(t)) - \int_t^1 
	 g(s)ds }{ \log(1/t)}+1$, and the hypothesis becomes $\limt 
	 \frac{\log(tg(t)) - \int_t^1 g(s)ds }{ \log(1/t)}$ $= 0$, that is 
	 for any $\eps>0$ there is $t_\eps>0$ s.t.  
	 \begin{equation}\label{e:liminf1}
		 \frac{\log(tg(t)) - \int_t^1 g(s)ds}{ \log(1/t)} \geq -\eps\; ,
	 \end{equation}
	 for all $t\in (0,t_\eps)$.  Observe now that $\frac{S(2t)}{ 
	 S(t)}=\e{\int_t^{2t} g(s)ds}$, so that 
	 $\liminf_{t\to0^+}\frac{S(2t)}{ S(t)}$ $=1$ is equivalent to $\limt 
	 \int_t^{2t} g(s)ds = 0$, and, as $0 \leq \int_t^{2t} g(s)ds \leq 
	 tg(t)$, our thesis will be proved as soon as we can show that 
	 $\limt tg(t)=0$.  So suppose, on the contrary, there are 
	 $t_0,c>0$ s.t. $tg(t)\geq c$, for all $t\in(0,t_0)$.  From 
	 inequality (\ref{e:liminf1}) it follows that
	 \begin{align*}
		 \frac{\log(tg(t))}{ \log(1/t)} & \geq \frac{\int_t^1 g(s)ds}{ 
		 \log(1/t)} - \eps \cr & \geq \frac{c\log(t_0/t)}{ \log(1/t)} 
		 + \frac{\int_{t_0}^1 g(s)ds}{ \log(1/t)} - \eps \cr & = c 
		 -\eps + \frac{c\log t_0 + \int_{t_0}^1 g(s)ds}{ \log(1/t)}
	 \end{align*}
	 for all $t\in (0,\min\{t_\eps,t_0\})$.  As 
	 $\lim_{t\to0^{+}}\frac{c\log t_0 + \int_{t_0}^1 g(s)ds}{ 
	 \log(1/t)} = 0$, we can choose $\eps<c/4$, and $t_1>0$ s.t. 
	 $\frac{\log(tg(t))}{ \log(1/t)}\geq \frac{c}{2}$, for all 
	 $t\in(0,t_1)$.  Then we have $g(t)\geq (1/t)^{c/2+1}$, for all 
	 $t\in(0,t_1)$, so that
	 \begin{equation*}
		 \frac{\int_t^1 g(s)ds}{ \log(1/t)} \geq \frac{\int_t^{t_1} 
		 (1/s)^{c/2+1}ds + \int_{t_1}^1 g(s)ds }{ \log(1/t)} 
	 \end{equation*}
	 from which it follows 
	 \begin{equation}\label{e:infinito}
		 \lim_{t\to0^{+}}\frac{\int_t^1 g(s)ds}{ \log(1/t)} = +\infty. 
	 \end{equation} 
	 As inequality (\ref{e:liminf1}) it equivalent to 
	 $\frac{\log(g(t)) - \int_t^1 g(s)ds}{ \log(1/t)} \geq 1 -\eps$, 
	 for all $t\in (0,t_\eps)$, it follows that $\frac{\log g(t)}{ 
	 \int_t^1 g(s)ds}\geq 1 + (1-\eps) \frac{\log(1/t) }{ \int_t^1 
	 g(s)ds} \to 1$, when $t\to0^+$.  Therefore, with $G(t):=\int_t^1 
	 g(s)ds$, we have $\frac{\log(-G'(t))}{ G(t)} \geq 1-\eps$, for 
	 all $t\in(0,t_2)$.  This implies $-G'(t)\e{-(1-\eps)G(t)}\geq 1$, 
	 which we integrate in $(0,t)$ for $t\leq t_2$, to obtain $t\leq 
	 \frac{1}{ 1-\eps}\ \e{-(1-\eps)G(t)}$ that is $G(t)\leq \frac{1}{ 
	 1-\eps} [\log(1/t) - \log(1-\eps)]$, and choosing $\eps$ small 
	 enough, we get $G(t)\leq 2\log(1/t)$, for all $t$ near $0$.  But 
	 this contrasts with equation (\ref{e:infinito}), so we are done.
 \end{proof}
 
 \begin{Thm}\label{inf-ecc} 
	 Let $a\in\mt$ be s.t. $\oinf{a}=1$.  Then $a$ is 
	 $\infty$-eccentric.
 \end{Thm}
 \begin{proof}
	 It is a consequence of the following Propositions.
 \end{proof}

 \begin{Prop}\label{1.9} 
	 Suppose $\m_a\not\in L^1(1,\infty)$ and $\liminf_{t\to\infty} 
	 \frac{\log\m_a(t)^{-1}}{\log t}=1$.  \\
	 Then $\liminf_{t\to\infty}\frac{S^\infty_a(2t)}{ 
	 S^\infty_a(t)}=1$.
 \end{Prop}
 \begin{proof}
	 Recall that $S(t)\equiv S^\infty_a(t)=\int_1^t \m_a(s)ds$, is 
	 positive, nondecreasing, concave, and $S(1)=0$, 
	 $\lim_{t\to\infty}S(t)=\infty$.  \\
	 Then $\log S$ is nondecreasing and concave, and 
	 $\lim_{t\to\infty}\log S=\infty$, so that $g:=\frac{S'}{ S}$ is 
	 positive, nonincreasing, $\lim_{t\to\infty}g(t)=0$, $\log S(t) = 
	 \int_1^t g(s)ds$, and $\int_1^\infty g(s)ds = 
	 \lim_{t\to\infty}\log S = \infty$.  Finally $\m_a=S'=gS$, so that 
	 $\frac{\log\m_a(t)^{-1}}{\log t} = -\frac{\log(tg(t))+\int_1^t 
	 g(s)ds }{ \log t}+1$, and the hypothesis becomes
	 \begin{equation}\label{e:liminf2}
		 \limsup_{t\to\infty} \frac{\log(tg(t)) + \int_1^t 
		 g(s)ds}{\log t} = 0
	 \end{equation}
	 Observe now that $\frac{S(2t)}{ S(t)}=\e{\int_t^{2t} g(s)ds}$, so 
	 that $\liminf_{t\to\infty}\frac{S(2t)}{ S(t)}=1$ is equivalent to 
	 $\liminf_{t\to\infty} \int_t^{2t} g(s)ds = 0$, and, as $0 \leq 
	 \int_t^{2t} g(s)ds \leq tg(t)$, our thesis will be proved as soon 
	 as we can show that $\liminf_{t\to\infty} tg(t)=0$.  So suppose, 
	 on the contrary, there are $t_0,c>0$ s.t. $tg(t)\geq c$, for all 
	 $t\in(t_0,\infty)$.  From equation (\ref{e:liminf2}) it follows 
	 that for all $\eps>0$ there is $t_\eps>0$ s.t. $\frac{\log(tg(t)) 
	 + \int_1^t g(s)ds}{ \log t} \leq \eps$, for all $t>t_\eps$, 
	 so that
	 \begin{align*}
		 \frac{\log(tg(t))}{ \log t} & \leq -\frac{\int_1^t 
		 g(s)ds}{ \log t} + \eps \cr & \leq -\frac{\int_1^{t_0} 
		 g(s)ds}{ \log t} - \frac{c\log(t/t_0)}{ \log t} + \eps \cr & 
		 = -c +\eps + \frac{c\log t_0 - \int_1^{t_0} g(s)ds}{ \log t}
	 \end{align*}
	 for all $t > \max\{t_\eps,t_0\}$.  As $\lim_{t\to\infty} 
	 \frac{c\log t_0 - \int_1^{t_0} g(s)ds}{ \log t} = 0$, we can 
	 choose $\eps<c/4$, and $t_1>\max\{t_\eps,t_0\}$ s.t. 
	 $\frac{\log(tg(t))}{ \log t}\leq -\frac{c}{2}$, for all $t>t_1$.  
	 Then we have $g(t)\leq t^{-c/2-1}$, for all $t>t_1$, so that $c/t 
	 \leq g(t) \leq t^{-c/2-1}$, for all $t>t_1$, and this is absurd.
 \end{proof}
 
 \begin{Prop}\label{1.10} 
	 Suppose $\m_a\in L^1(1,\infty)$ and $\liminf_{t\to\infty} 
	 \frac{\log\m_a(t)^{-1}}{\log t}=1$.  \\
	 Then $\limsup_{t\to\infty}\frac{S^\infty_a(2t)}{ 
	 S^\infty_a(t)}=1$.
 \end{Prop}
 \begin{proof}
	 Setting $S:= S^\infty_a$, as $S$ is a decreasing and convex 
	 function, $1 \geq \frac{S(2t)}{S(t)} \geq 1 - 
	 \frac{t\m(t)}{S(t)}$, so that 
	 $\liminf_{t\to\infty}\frac{t\m(t)}{S(t)}=0$ implies the thesis.  
	 Let us assume on the contrary that there are $c, t_{0}>0$ s.t. 
	 $\frac{t\m(t)}{S(t)}\geq c$, for all $t>t_{0}$.  Then 
	 $\frac{S'(t)}{S(t)} = -\frac{\m(t)}{S(t)} \leq -\frac{c}{t}$, so 
	 that $S(t) \leq S(t_{0})\left(\frac{t}{t_{0}}\right)^{-c}= 
	 kt^{-c}$, for some positive $k$.  Then $t\m(2t) \leq 
	 \int_{t}^{2t} \m(s)ds \leq \int_{t}^{\infty} \m(s)ds \leq 
	 kt^{-c}$, for all $t\geq t_{0}$, therefore $\m(t) \leq 
	 k't^{-c-1}$, for all $t\geq 2t_{0}$, so that 
	 $\liminf_{t\to\infty} \frac{\log \m(t)}{\log(1/t)} \geq 1+c$ 
	 contrary to the hypothesis.
 \end{proof}
 
 \medskip

 Because of their importance in determining the eccentricity 
 properties of an operator $a$, let us compute the orders of $a$ in 
 terms of the asymptotics of its distribution function $\l_{a}$.

 \begin{Prop}\label{1.8} 
	 $\ozero{a}$ and $\oinf{a}$ are also given by
	 \begin{align*}
		 \ozero{a} & = 
		 \left(\limsup_{s\to\infty}\frac{\log\l_a(s)}{\log(1/s)}\right)^{-1} \; , 
		 \\
		 \oinf{a} & = 
		 \left(\limsup_{s\to0^{+}}\frac{\log\l_a(s)}{\log(1/s)}\right)^{-1} 
		 \; .
	 \end{align*}
 \end{Prop}
 \begin{proof}
	 Let us set $\m:=\m_a$, $\l:=\l_a$, and assume $\m$ is never 0 and 
	 $\lim_{t\to0}\m(t)=+\infty$, otherwise the proof is obvious.\\
	 Let $G_1=\left\{(x,y)\in \br^{2}_{+} : y=\m(x) \right\}$, 
	 $G_2=\left\{(x,y)\in \br^{2}_{+} : x=\l(y) \right\}$.  Then
	 \begin{equation*}
		 \liminf_{t\to\infty} \frac{\log\m(t)}{\log(1/t)} = 
		 \lim_{t\to\infty} \inf_{\substack{(x,y)\in G_1 \\
		 x>t}}  \frac{\log y}{\log(1/x)}  = \lim_{s\to0} 
		 \inf_{\substack{(x,y)\in G_1\\ y<s}}  \frac{\log 
		 y}{\log(1/x)} 
	 \end{equation*}
	 and
	 \begin{equation*}
		 \liminf_{s\to0} \frac{\log s}{\log(1/\l(s))} = \lim_{s\to0} 
		 \inf_{\substack{(x,y)\in G_2 \\ y<s}}  \frac{\log 
		 y}{\log(1/x)}  = \lim_{t\to\infty} 
		 \inf_{\substack{(x,y)\in G_2 \\
		 x>t}}  \frac{\log y}{\log(1/x)} .
	 \end{equation*}
	 Also 
	 $$
	 \lim_{t\to\infty} \inf_{\substack{(x,y)\in G_1 \\
	 x>t}}  \frac{\log y}{\log(1/x)}  = 
	 \lim_{t\to\infty} \inf_{\substack{(x,y)\in\ov{G_1} \\ 
	 x>t}} \frac{\log y}{\log(1/x)},
	 $$
	 $$
	 \lim_{s\to0} \inf_{\substack{(x,y)\in G_1 \\ 
	 y<s}}  \frac{\log y}{\log(1/x)}= \lim_{s\to0} 
	 \inf_{\substack{(x,y)\in\ov{G_1} \\ y<s}} \frac{\log 
	 y}{\log(1/x)},
	 $$
	 $$
	 \lim_{t\to\infty}\inf_{\substack{(x,y)\in G_2 \\ 
	 x>t}} \frac{\log y}{\log(1/x)}= 
	 \lim_{t\to\infty}\inf_{\substack{(x,y)\in\ov{G_2} \\ x>t}} \frac{\log 
	 y}{\log(1/x)},
	 $$
	 and
	 $$
	 \lim_{s\to0} \inf_{\substack{(x,y)\in G_2 \\ 
	 y<s}}  \frac{\log y}{\log(1/x)} = \lim_{s\to0} 
	 \inf_{\substack{(x,y)\in\ov{G_2} \\ y<s}} \frac{\log 
	 y}{\log(1/x)}.
	 $$
	 Define, for any $y>0$, $[\ell_y,r_y] := \left\{x>0 : 
	 (x,y)\in\ov{G_1}\right\}$, and, for any $x>0$, 
	 $[d_x,u_x]=\left\{y>0 : (x,y)\in\ov{G_2}\right\}$.  Then 
	 $(x,y)\in\ov{G_1}\cap\ov{G_2}$ $\iff$ $(x,y)\in\ov{G_1}$ and 
	 $x\in\left\{\ell_y,r_y\right\}$ $\iff$ $(x,y)\in\ov{G_2}$ and 
	 $y\in\left\{d_x,u_x\right\}$.  Indeed, defining $\l^{-}(y) := 
	 \lim_{s\to y^{-}} \l(s)$, and $\m^{-}(x) := \lim_{t\to x^{-}} 
	 \m(t)$, one proves that $\ov{G_{1}(\m^{-})} = \ov{G_{1}(\m)}$ and 
	 $\ov{G_{2}(\l^{-})} = \ov{G_{2}(\l)}$, and the rest follows 
	 easily. \\
	 Finally, since 
	 $\frac{\log\m(t)}{\log(1/r_t)}\leq\frac{\log\m(t)}{\log(1/t)}$ 
	 one gets, taking the $\liminf$,
	 $$
	 \lim_{t\to\infty}\inf_{\substack{(x,y)\in\ov{G_1} \\ x>t}}\frac{\log 
	 y}{\log(1/x)}= 
	 \lim_{t\to\infty}\inf_{\substack{(x,y)\in\ov{G_1}\cap\ov{G_2} \\ x>t}} \frac{\log 
	 y}{\log(1/x)}.
	 $$
	 Analogously, since $\frac{\log u_s}{\log(1/\l(s))}\leq\frac{\log 
	 s}{\log(1/\l(s))}$ one gets, taking the $\liminf$,
	 $$
	 \lim_{s\to0}\inf_{\substack{(x,y)\in\ov{G_2} \\ y<s}}\frac{\log 
	 y}{\log(1/x)}= 
	 \lim_{s\to0}\inf_{\substack{(x,y)\in\ov{G_1}\cap\ov{G_2} \\ y<s}}\frac{\log 
	 y}{\log(1/x)}.
	 $$
	 Putting all this together one gets
	 \begin{align*}
		 \liminf_{t\to\infty}\frac{\log\m(t)}{\log(1/t)} & = 
		 \liminf_{s\to0}\frac{\log s}{\log(1/\l(s))}= 
		 \left(\limsup_{s\to0}\frac{\log(1/\l(s))}{\log s}\right)^{-1} 
		 \\
		 & = \left(\limsup_{s\to0}\frac{\log\l(s)}{\log(1/s)}\right)^{-1}.
	 \end{align*}
	 The equality
	 $$
	 \liminf_{t\to0}\frac{\log\m(t)}{\log(1/t)}= 
	 \left(\limsup_{s\to\infty}\frac{\log\l(s)}{\log(1/s)}\right)^{-1}
	 $$
	 is proved in the same way (using $\ell_{t}$ and $d_{s}$).
 \end{proof}

\section{Some results on noncommutative geometric measure theory}\label{sec2}

 In this section we shall discuss some definitions of dimension in 
 noncommutative geometry in the spirit of geometric measure theory.

 As it is known, the measure for a noncommutative manifold is defined 
 via a singular trace applied to a suitable power of some geometric 
 operator (e.g. the Dirac operator of the spectral triple of Alain 
 Connes).  Connes showed that such procedure recovers the usual volume 
 in the case of compact Riemannian manifolds, and more generally the 
 Hausdorff measure in some interesting examples \cite{Co}.

 Let us recall that $(\ca,D,\ch)$ is called a spectral triple when 
 $\ca$ is an algebra acting on the Hilbert space $\ch$, $D$ is a self 
 adjoint operator on the same Hilbert space such that $[D,a]$ is 
 bounded for any $a\in\ca$, and $D$ has compact resolvent.  In the 
 following we shall assume that $0$ is not an eigenvalue of $D$, the 
 general case being recovered by replacing $D$ with 
 $D|_{\ker(D)^\perp}$.  Such a triple is called $d^+$-summable, $d\in 
 (0,\infty)$, when 
 $|D|^{-d}$ belongs to the ideal
 \begin{equation}\label{macaev}
	 \cl^{1,\infty}:=\{a\in\ck(\ch):\sum_{k=1}^{n}\m_k(a)=O(\log n)\},
 \end{equation}
 where, in case of compact operators, we denote the 
 non-increasing rearrangement of $a$ by $\m_{k}(a)$, instead of 
 $\m_{a}(k)$, to conform with tradition. \\ 
 The noncommutative version of the integral on functions is given 
 by the formula $\t_\omega(a|D|^{-d})$, where $\t_\omega$ is a 
 Dixmier trace, i.e. a singular trace summing logarithmic divergences.  
 Of course the preceding formula does not guarantee the non-triviality 
 of the integral, and in fact cohomological assumptions in this 
 direction have been considered \cite{Co}.  We are interested in 
 different conditions for non-triviality.  In this connection, 
 introducing the space $\cl^{1,\infty}_0$, where the $O(\log n)$ in 
 (\ref{macaev}) is replaced by $o(\log n)$, we observe that the 
 previous noncommutative integration is always trivial when $|D|^{-d}$ 
 belongs to $\cl^{1,\infty}_0$.

 \begin{Lemma}\label{Lemma:dim} 
	 Let $(\ca,D,\ch)$ be a spectral triple. Then
	 $$
	 \inf\{d>0:|D|^{-d}\in\cl^{1,\infty}_0\}= 
	 \sup\{d>0:|D|^{-d}\not\in\cl^{1,\infty}\}.
	 $$
 \end{Lemma}
 \begin{proof}
	 Let $d^+=\inf\{d>0:|D|^{-d}\in\cl^{1,\infty}_0\}$, $d<d^+$.  Then 
	 $|D|^{-d'}\not\in\cl^{1,\infty}_0$ for $d'=\frac{d+d^+}2$, i.e. 
	 there exists a subsequence $n_k$ such that
	 $$
	 \lim_{k\to\infty}\frac1{\log n_k}\sum_{j=1}^{n_k}\m_j(|D|^{-d'})=\ell>0.
	 $$
	 Then, setting $\eps=1-\frac{d}{d'}>0$, for any $m\in\bn$ we have
	 \begin{align*}
		 \lim_{k\to\infty}\frac1{\log n_k}\sum_{j=1}^{n_k}\m_j(|D|^{-d}) 
		 &=\lim_{k\to\infty}\frac1{\log 
		 n_k}\sum_{j=m}^{n_k}\m_j(|D|^{-d})\cr 
		 &=\lim_{k\to\infty}\frac1{\log n_k} 
		 \sum_{j=m}^{n_k}\frac{\m_j(|D|^{-d'})}{\m_j(|D|^{-d'})^\eps}\cr 
		 &\geq \frac1{\m_{m}(|D|^{-d'})^\eps} \lim_{k\to\infty}\frac1{\log 
		 n_k}\sum_{j=m}^{n_k}\m_j(|D|^{-d'})\cr 
		 &=\frac{\ell}{\m_{m}(|D|^{-d'})^\eps}.
	 \end{align*}
	 By the arbitrariness of $m$, the limit is $\infty$, i.e. 
	 $|D|^{-d}\not\in\cl^{1,\infty}$, which implies the thesis.
 \end{proof}

 These results, together with the examples by Connes and Sullivan 
 \cite{Co}, justify the following definition.

 \begin{Dfn} 
	 Let $(\ca,D,\ch)$ be a spectral triple.  We shall call the 
	 functional $a\mapsto\t_\omega(a |D|^{-\a})$ the $\a$-dimensional 
	 Hausdorff measure, and the number
	 $$
	 d_H(\ca,D,\ch)=\inf\{d>0:|D|^{-d}\in\cl^{1,\infty}_0\}= 
	 \sup\{d>0:|D|^{-d}\not\in\cl^{1,\infty}\}
	 $$
	 the Hausdorff dimension of the spectral triple.
 \end{Dfn}

 Let us observe that the $d$-dimensional Hausdorff measure depends on 
 the generalized limit procedure $\omega$, however all such 
 functionals coincide on measurable operators in the sense of Connes 
 \cite{Co}.  As in the commutative 
 case, the Hausdorff dimension is the supremum of the $d$'s such that 
 the $d$-dimensional Hausdorff measure is everywhere infinite and the 
 infimum of the $d$'s such that the $d$-dimensional Hausdorff measure 
 is identically zero.

 Concerning the non-triviality of the $d_H$-dimensional Hausdorff 
 measure, we have the same situation as in the classical case.

 \begin{Prop} 
	 Let $(\ca,D,\ch)$ be a spectral triple with finite non-zero 
	 Hausdorff dimension $d_H$.  Then the $d_H$-dimensional Hausdorff 
	 measure is the only possibly non-trivial functional on $\ca$ 
	 among the Hausdorff measures.
 \end{Prop}
 \begin{proof} 
	 The result obviously follows by Lemma~\ref{Lemma:dim} and the 
	 definition of the Hausdorff measures.
 \end{proof}

 According to the previous result, a non-trivial Hausdorff measure is 
 unique but does not necessarily exist.  In fact, if the eigenvalue 
 asymptotics of $D$ is e.g. $n\log n$, the Hausdorff dimension is one, 
 but the 1-dimensional Hausdorff measure gives the null functional.

 We shall now propose another spectral dimension, for which the 
 situation is somewhat the opposite.  If we consider all singular 
 traces, not only the logarithmic ones, and the corresponding 
 functionals on $\ca$, we shall show that there exists a non trivial 
 functional associated with such a dimension, but such property does 
 not characterize this dimension.

 \begin{Dfn} 
	 Let $(\ca,D,\ch)$ be a spectral triple. We shall call the number
	 \begin{equation*}
		 d_{B}(\ca,D,\ch)=\oinf{D^{-1}}^{-1} 
		 =\left(\liminf_{n\to\infty}\frac{\log \m_n(D)}{\log n 
		 }\right)^{-1}
	 \end{equation*}
	 the box dimension of the spectral triple.
 \end{Dfn}

 \begin{Prop} 
	 Let $(\ca,D,\ch)$ be a spectral triple with finite non-zero box 
	 dimension $d$.  Then $|D|^{-d}$ is singularly traceable, namely it 
	 gives rise to a singular trace $\t$ which is non-trivial on the 
	 ideal generated by $|D|^{-d}$.  In particular the functional 
	 $a\mapsto\t(a |D|^{-d})$ is a non-trivial trace state on the 
	 algebra $\ca$.
 \end{Prop}
 \begin{proof} 
	 By definition $\oinf{D^{-1}}=d^{-1}$ hence (cf.  Remark 
	 \ref{rem1.6}) $\oinf{|D|^{-d}}=1$, therefore, by Theorem 
	 \ref{inf-ecc}, $|D|^{-d}$ is eccentric and finally, by Theorem 
	 \ref{1.4}, we get the existence of a singular trace $\t$.  The 
	 trace property for the functional $a\mapsto\t(a |D|^{-d})$ is 
	 proved as in \cite{CiGS1}. 
 \end{proof}

 \begin{rem}\label{uniqueness} 
	 We call the number $d_B$ a dimension since it is related to
	 the existence of a non-trivial geometric measure.
	 Proposition~\ref{Prop:unique} shows that under suitable
	 regularity conditions of the eigenvalue sequence $\mu_n(D)$
	 such request determines $d_B$ uniquely, and $d_B$ coincides
	 with $d_H$.  However this is not true in general.  In a
	 following example we describe some selfadjoint operators $D$ for
	 which the numbers $d_H$ and $d_B$ are different and both
	 $|D|^{-d_H}$ and $|D|^{-d_B}$ are singularly traceable.
 \end{rem}

 \begin{Prop}\label{Prop:unique} 
	 Let $(\ca,D,\ch)$ be a spectral triple with finite non-zero box 
	 dimension. \\
	 \item{$(a)$} If there exists $\lim\frac{\log\m_n(D^{-1})}{\log 
	 1/n}$, then $d_B=d_H$. \\
	 \item{$(b)$} If there exists 
	 $\lim\frac{\m_n(D^{-1})}{\m_{2n}(D^{-1})}$, $d_B$ is 
	 characterized by the property that $|D|^{-d_{B}}$ is singularly 
	 traceable, and $d_B=d_H$.
 \end{Prop}
 \begin{proof} 
	 $(a)$.\quad As $d_{H}(|D|^{\a}) = \frac{1}{\a} d_{H}(|D|)$, and 
	 $d_{B}(|D|^{\a}) = \frac{1}{\a} d_{B}(|D|)$, for any $\a>0$, we 
	 may restrict to the case $d_B=1$.  By hypothesis we have that for 
	 any $\eps>0$
	 $$
	 \left(\frac1n\right)^{1+\eps} \leq\m_n(D^{-1}) 
	 \leq\left(\frac1n\right)^{1-\eps}.
	 $$
	 for sufficiently large $n$.  As a consequence, if $\l>1$,
	 $$
	 \m_n(|D|^{-\l})\leq\left(\frac1n\right)^{\frac{\l+1}2}
	 $$
	 hence it is a summable sequence, which implies 
	 $\frac{\sum_{k=1}^n\m_n(|D|^{-\l})}{\log n}=0$, i.e. $d_H\leq1$.  
	 \\
	 Conversely, if $\l<1$, 
	 $$
	 \m_n(|D|^{-\l})\geq \left(\frac1n\right)^{\frac{\l+1}2}
	 $$
	 and this implies that $\frac{\sum_{k=1}^n\m_n(|D|^{-\l})}{\log 
	 n}=\infty$, i.e. $d_H\geq1$.  The thesis follows.  \\
	 $(b)$.\quad For any $n\in\bn$, let $k_n\in\bn$ be such that 
	 $2^{k_n}\leq n<2^{k_n+1}$, and write $\m_n$ for $\m_n(D^{-1})$ and 
	 assume for simplicity that $\m_1=1$.  Then
	 \begin{align*}
		 -\frac1{k_n\log2}\sum_{j=0}^{k_n}\log\frac{\m_{2^{j+1}}}{\m_{2^j}} 
		 &=-\frac{\log\m_{2^{k_n+1}}}{k_n\log2}\cr &\geq\frac{\log 
		 \m_n}{\log 1/n}\cr 
		 &\geq-\frac{\log\m_{2^{k_n}}}{(k_n+1)\log2}\cr 
		 &=-\frac1{(k_n+1)\log2}\sum_{j=0}^{k_n-1}\log\frac{\m_{2^{j+1}}}{\m_{2^j}}.
	 \end{align*}
	 Taking the limit for $n\to\infty$ one gets that 
	 $\lim\frac{\log\m_n(D^{-1})}{\log 1/n}$ exists, hence $d_B=d_H$ 
	 by $(a)$, and also
	 $$
	 d_B^{-1}=-\frac{\log\left(\lim_k\frac{\m_{2k}}{\m_k}\right)}{\log2},
	 $$
	 namely
	 $$
	 \lim_n\frac{\m_{2n}(D^{-1})}{\m_n(D^{-1})}=2^{-1/d_B}.
	 $$
	 Assume for the moment that $\m_{n}\not\in \ell^{1}$, and denote 
	 by $s_{n}(|D|^{-d}) := \sum_{k=1}^{n}\m_k(|D|^{-d})$ (the same as 
	 $S^{\infty}_{a}(n)$ of section 1).  Then, by a 
	 Cesaro theorem,
	 \begin{align*}
		 \lim_n\frac{s_{2n}(|D|^{-d})}{s_n(|D|^{-d})} 
		 &=\lim_n\frac{\sum_{k=1}^{2n}\m_k(|D|^{-d})}{\sum_{k=1}^{n}\m_k(|D|^{-d})} 
		 =\lim_n\frac{2\m_{2n}(|D|^{-d})}{\m_n(|D|^{-d})}\cr 
		 &=\lim_n\frac{2\m_{2n}(|D|^{-1})^d}{\m_n(|D|^{-1})^d} 
		 =2\left(\lim_n\frac{\m_{2n}(|D|^{-1})}{\m_n(|D|^{-1})}\right)^d 
		 =2^{1-d/d_B}.
	 \end{align*}
	  Therefore $|D|^{-d}$ is eccentric if and only if such limit is 
	  one, i.e. when $d=d_B$.  If $\m_{n}\in 
	  \ell^{1}$, then denoting by $s_{n}(|D|^{-d}) := 
	  \sum_{k=n}^{\infty}\m_k(|D|^{-d})$ (the same as 
	  $S^{\infty}_{a}(n)$ of section 1), the calculation above, 
	  suitably modified, shows that $|D|^{-d}$ is eccentric if and 
	  only $d=d_B$.
 \end{proof}

 \begin{exmp}\label{ex:dim}
	 Let us construct a family of Dirac operators $D_{\l}$, $\l>1$, s.t. 
	 $d_{B}(D_{\l})=1$, $d_{H}(D_{\l})=\l$, and the $\l$-dimensional 
	 Hausdorff measure is non-trivial.
	 Since the dimensions and the singular traceability property depend 
	 only on the eigenvalue sequence $\m_n(\l):=\m_n(|D_{\l}|^{-1})$, we shall 
	 concentrate only on the construction of the sequence $\m_n(\l)$.\\
	 Let $a_k$ be any increasing diverging sequence, $a_1=0$, and set 
	 $\m_n=\e{-a_k}$ when $\e{a_k}\leq n<\e{a_{k+1}}$.  Then
	 $$
	 d_B^{-1}=\liminf_n\frac{\log \m_n}{\log 1/n}= \lim_k \frac{\log 
	 \m_{[\e{a_k}]+1}}{\log 1/([\e{a_k}])}=1
	 $$
	 where $[\cdot]$ denotes the integer part.  If, setting 
	 $\s_{n,\l}=\s_n(|D_{\l}|^{-\l})=\sum_{k=1}^n \m_k(\l)^{\l}$, we show that
	 $$
	 \limsup_n\frac{\s_{n,\l}}{\log n}
	 $$
	 is finite non-zero, this shows at once that $d_H(D_{\l})=\l$ and that 
	 there exists a non trivial logarithmic singular trace on 
	 $|D_{\l}|^{-\l}$ \cite{Varga,AGPS}. \\
	 Now, for any $\l>1$, set $a_k:=\l^k-\frac{\log\l}{\l-1}k$ and 
	 observe that, with this choice, $a_{j+1}-\l a_j=j\log \l - 
	 \log\l/(\l-1)$.  Then
	 \begin{align*}
		 \limsup_n \frac {\s_{n,\l}} {\log n} & = 
		 \lim_k \frac {\s_{[\e{a_k}],\l}} {\log [\e{a_k}]} = 
		 \lim_k \frac {\s_{[\e{a_k}],\l}} {a_k} \cr
		 & = \lim_k 
		 a_k^{-1}\sum_{j=1}^{k-1}(\e{a_{j+1}}-\e{a_j})\e{-\l a_j} \cr
		 & = \lim_k \l^{-k}\sum_{j=1}^{k-1}\e{a_{j+1}-\l a_j} = \lim_k 
		 \l^{-k}\sum_{j=1}^{k-1}\l^j\l^{\frac{1}{1-\l}} \cr
		 & = \lim_k \l^{-k}\l^{\frac{1}{1-\l}}\frac{\l^k-1}{\l-1} 
		 =\frac {\l^{\frac{1}{1-\l}}} {\l-1}.
	 \end{align*}
 \end{exmp}

 \begin{rem} 
	 Example~\ref{ex:dim} describes situations where the two geometric 
	 spectral dimensions considered here are different, and give rise 
	 to different (non trivial) geometric integrations. \\
	 For the spectral triples whose Dirac operator has a spectral 
	 asymptotics like $n^{\a}(\log n)^{\b}$ instead, we have 
	 $d_B=d_H=1/\a$, namely the two dimensions coincide, and the 
	 uniqueness result of the preceding Proposition applies.  However, 
	 the nontrivial singular trace associated with $|D|^{-d_B}$ by 
	 Theorem \ref{1.4} is a logarithmic trace if and only if $\b=1$.  
	 In this sense, the singular traces associated with a generic 
	 eccentric operator generalize the logarithmic trace in the same 
	 way in which the Besicovitch measure theory generalizes the 
	 Hausdorff measure theory.
 \end{rem}


\section{Novikov-Shubin invariants as asymptotic dimensions}

 In this section we apply the theory of polynomial orders of 
 operators introduced in section \ref{sec:singular} to a geometric 
 operator (the Laplacian on $k$-forms) and show that these orders 
 give topological information on the manifold.

 \subsection{Weyl's asymptotics via Atiyah's trace on covering 
 ma\-ni\-folds}
 \label{subsec:Atiyah}
 
 Let $M$ be a complete connected Riemannian $n$-manifold, and $G$ an 
 infinite discrete group of isometries of $M$.  Suppose that $G$ acts 
 freely ($i.e.$ any $g\in G$, $g\neq e$, acts without fixed points), 
 properly discontinuously, and that $X:=M/ G$ is a compact manifold.  
 Let $\cf$ be a fundamental domain for $G$, that is (\cite{Atiyah}, 
 page 52) an open subset of $M$, disjoint from all its translates by 
 $G$, and such that $M\setminus \cup_{g\in G} g(\cf)$ has measure 
 zero.  Let $L^2\La^k(M)$ be the Hilbert space of square-integrable 
 $k$-forms on $M$, w.r.t. the volume measure, then $L^2\La^k(M)\cong 
 \ell^2(G)\otimes L^2\La^k(X)$.  $G$ acts on $L^2\La^k(M)$ as left 
 translation operators $(L_gu)(x) := u(g^{-1}x)$.  Let $\cam_k\equiv 
 \cam_k(M,G)$ be the von Neumann algebra of bounded $G$-invariant 
 operators, so that $\cam_k\cong \car(G)\otimes \cb(L^2\La^k(X))$ and 
 $\cam_k' = \{L_g: g\in G\}''\cong \cl(G)\otimes\bc$, where $\car(G),\ 
 \cl(G)$ are the right, resp.  left, regular representations of $G$.  
 Any self-adjoint $G$-invariant operator on $L^2\La^k(M)$ is 
 affiliated with $\cam_k$.

 By the previous isomorphism, $\cam_k$ inherits a trace $Tr_G= 
 \tau_G\otimes Tr$, and we quote a result in \cite{Atiyah}, which 
 gives a more explicit description of $Tr_{G}$.  

 \begin{Prop}\label{2.3}  
	 Let $A\in\cam_k$ be a positive self-adjoint operator, with a 
	 $C^\infty$ kernel $A(x,y)$.  Then $A\in \cl^1(\cam_k,Tr_G)$, and 
	 $Tr_G(A)=\int_\cf tr A(x,x) dvol(x)$, where $tr$ is the usual 
	 matrix trace.
 \end{Prop}

 The Laplacian $\D_k$ acting on exterior $k$-forms on $M$ is 
 essentially self-adjoint as an operator on $L^2\La^k(M)$ 
 \cite{Chernoff}, and we use the same notation for its closure.  Let 
 $\D_k = \int tdE_k(t)$, be its spectral decomposition; then 
 $e_k(t,\cdot,\cdot)$, the Schwartz kernel of $E_k(t)$, belongs to 
 $C^\infty(M\times M)$, and we have for the spectral distribution 
 function $N_k(t) := Tr_G(E_k(t)) = \int_\cf tr e_k(t,x,x)dvol(x)$.  
 $N_k$ is an increasing function on $\br$ which vanishes on 
 $(-\infty,0)$.  \par

 
 We are now in a position to make explicit the topological information 
 contained in $\oinf{\D_{k}^{-1}}$. We need a lemma.

 \begin{Lemma}\label{2.4} 
	 $\l_{\D_k^{-1}}(t)=N_k(1/t)-b_k$.
 \end{Lemma}
 \begin{proof}
	 \begin{align*}
		 \l_{\D_k^{-1}}(t) & = Tr_G(E_{(t,\infty)}(\D_k^{-1})) = 
		 Tr_G(\chi_{(t,\infty)}(\D_k^{-1})) \cr & = 
		 Tr_G(\chi_{(0,1/t)}(\D_k)) = Tr_G(E_{(0,1/t)}(\D_k)) = 
		 N_k(1/t)-b_k.
	 \end{align*}
 \end{proof}

 \begin{Thm} \label{t:inf-ecc} 
	 $dim(M) = 2(\oinf{\D_{k}^{-1}})^{-1}$.  As a consequence 
	 $\D_{k}^{-n/2}$ is $\infty$-eccentric, and gives rise to a 
	 singular trace on $\cam_{k}$.
 \end{Thm} 
 \begin{proof}
	 Recall from \cite{ES}, equation (4.5), that $N_{k}(t)\sim \b 
	 t^{n/2}$, as $t\to\infty$, where $\b\neq0$.  Then from Lemma 
	 \ref{2.4} it follows $\l_{\D_{k}^{-1}}(s) \sim \b s^{-n/2}$, as 
	 $s\to0$, so that the thesis follows from Proposition \ref{1.8}.
 \end{proof}

 \subsection{Novikov-Shubin invariants as asymptotic dimensions}
 \label{subsec:NSinvariants}
 
 Novikov and Shubin \cite{NS1}, \cite{NS2} have studied the asymptotic 
 behaviour of $N_k(t)$ as $t\to0$, which, through the efforts of 
 Efremov-Shubin \cite{ES}, Lott \cite{Lott}, and Gromov-Shubin 
 \cite{GS}, has been proved to be a homotopy invariant.  We want to 
 show that the Novikov-Shubin invariants are asymptotic dimensions, so 
 we need some notation.

 Let $\th_k(t) := Tr_G(\e{-t\D_k}) = \int \e{-st}dN_k(s)$ be the 
 Laplace-Stieltjes transform of $N_k(t)$.  $\th_k$ is a decreasing 
 positive function on $(0,\infty)$.

 Gromov and Shubin introduced (weak) Novikov-Shubin numbers, which, using 
 Lott normalization \cite{Lott}, are defined as follows
 \begin{itemize}
	 \itm{i} $\underline{\a}_k \equiv \underline{\a}_k(M,G) := 2\liminf_{t\to0} 
	 \frac{\log(N_k(t)-b_k) }{ \log t} = 2\liminf_{t\to\infty} 
	 \frac{-\log(\th_k(t)-b_k) }{ \log t}$ 
	 \itm{ii} $\a_k \equiv \a_{k}(M,G)
	 := 2\limsup_{t\to0} \frac{\log(N_k(t)-b_k) }{ \log t}$ \itm{iii} 
	 $\a'_k \equiv \a'_{k}(M,G) := 2\limsup_{t\to\infty} 
	 \frac{-\log(\th_k(t)-b_k) }{ \log t}$
 \end{itemize}
 where $b_k := \lim_{t\to0} N_k(t)$ are the so-called $L^2$-Betti 
 numbers, and are homotopy invariant \cite{Dodziuk}.  Gromov and 
 Shubin showed that these numbers are $G$-homotopy invariants of $M$.  
 \\
 In analogy with the definition given in \cite{GI2}, we call asymptotic
 spectral dimension of the covering manifold $M$ with structure group
 $G$ the number
 $$
 d_{\infty}(M,G,\D_{k}) := 2(\ozero{\D_{k}^{-1}})^{-1}.
 $$
 Then
 
 \begin{Thm}\label{3.2} 
	 Let $k$ be s.t. $0<\a_k<\infty$.  Then 
	 $\a_k = d_{\infty}(M,G,\D_{k})$.  Therefore 
	 $\D_k^{-\a_k/2}$ is $0$-eccentric, and gives rise to a 
	 non-trivial singular trace on $\cam_{k}$.
 \end{Thm}
 \begin{proof}
	 From Proposition \ref{1.8} and Lemma \ref{2.4} it follows that
	 $$
	 d_{\infty}(M,G,\D_{k}) = 2 \ozero{\D_{k}^{-1}}^{-1}=
	 2 \limsup_{t\to0}\frac{\log(N_k(t)-b_k)}{\log t} = \a_k.
	 $$
	 Therefore, if $0<\a_k<\infty$, 
	 $\ozero{\D_{k}^{-\a_k/2}}=1$
	 and the thesis follows from Theorem \ref{1.5}.
 \end{proof}
 
 \begin{rem} If $k=0$, a result by Varopoulos \cite{Varopoulos1} shows
 	that $\a_{0}(M,G) = {\rm growth\;} G$. Since the growth of $G$
 	coincides with the asymptotic metric dimension of $M$
 	\cite{GI2}, we obtain that the $0$-th Novikov-Shubin number
 	coincides with the asymptotic metric dimension.
 \end{rem}

 \subsection{Relation between Novikov-Shubin invariants and the 
 asymptotic dimension of the heat semigroups}\label{subsec:semigroupformula}

 Based on the notion of dimension at infinity due to Varopoulos, 
 Saloff-Coste, Coulhon \cite{VSC}, see also \cite{Cou}, we define the 
 asymptotic dimension of a semigroup of bounded operators on a measure 
 space.
 
 Let $(X,\cam,\mu)$ be a measure space, $V$ a finite dimensional (real 
 or complex) vector space, $L^{p}(X,\cam,\mu;V)$ the Lebesgue space 
 of $V$-valued functions. 

 \begin{Dfn}\label{2.2.1} 
	 Let $T_t : L^1(X,\cam,\m;V)\to L^\infty(X,\cam,\m;V)$ be a 
	 semigroup of bounded operators.  Then we set
	 $$
	 \ad{T} := \liminf_{t\to\infty} \frac{2\log \|T_t\|_{1\to\infty} 
	 }{ \log \left(\frac{1}{t}\right)} .
	 $$ 
 \end{Dfn}

 \begin{Thm}\label{2.2.2} {\rm (\cite{VSC}, Theorem II.4.3)} \\
	 Let $T_t\in \cb(L^1(X,\cam,\m;V)\cap L^\infty(X,\cam,\m;V))$ and 
	 assume it extends to a semigroup on $L^p$, for any 
	 $p\in[1,\infty]$, of class $C^0$ if $p<\infty$.  Suppose moreover 
	 that $T_t$ is equicontinuous on $L^1$ and $L^\infty$, bounded 
	 analytic on $L^2$, and $\|T_1\|_{1\to\infty}<\infty$.  Denote by 
	 $A$ the generator of the semigroup, and by $\cd:= {\rm span\ } \{ 
	 \int_0^\infty \f(t)T_tf dt\ :\ \f\in C^\infty_0(0,\infty),\ 
	 f\in\ L^\infty(X,\cam,\m;V),\ \m\{f\neq0\}<\infty \}$.  Then for 
	 any $n>0$, and $0<\a<\frac{n}{2}$, the following are equivalent 
	 \itm{i} $\|f\|_{2n/(n-2\a)} \leq C(\|A^{\a/2}f\|_2 + 
	 \|A^{\a/2}f\|_{2n/(n-2\a)})$, $f\in \cd$ 
	 \itm{ii} $\|T_1f\|_{2n/(n-2\a)} \leq C \|A^{\a/2}f\|_2$, $f\in \cd$  
	 \itm{iii} $\|T_t\|_{1\to\infty} \leq Ct^{-n/2}$, $t\in[1,\infty)$.
 \end{Thm}

 \begin{Prop}\label{2.2.3} 
	 Let $T_t\in \cb(L^1(X,\cam,\m;V)\cap L^\infty(X,\cam,\m;V))$ and 
	 assume it extends to a semigroup on $L^1$ of class $C^{0}$, and 
	 that $\|T_{1}\|_{1\to\infty}<\infty$.  Then the following 
	 are equivalent
	 \itm{i} $\|T_t\|_{1\to\infty} \leq Ct^{-n/2}$, $t\geq1$
	 \itm{ii} $\|T_t\|_{1\to\infty} \leq Ct^{-n/2}$, $t\geq t_0>1$. 
 \end{Prop}
 \begin{proof}
	 $(ii)\imply(i)$ Let $t\in(1,t_{0}]$ and observe that 
	 $\|T_t\|_{1\to\infty} = \|T_1T_{t-1}\|_{1\to\infty} \leq 
	 \|T_1\|_{1\to\infty}\|T_{t-1}\|_{1\to1} \leq 
	 k\|T_1\|_{1\to\infty} =: M$, where $k:= \sup_{t\in [0,t_{0}]} 
	 \|T_t\|_{1\to1}<\infty$ because $T_t$ is a semigroup of class 
	 $C^{0}$ on $L^1$.  So that, with $C_0:= \max\{ C, Mt_0^{n/2} \}$, 
	 we get the thesis.  
 \end{proof}\medskip

 \begin{Prop}\label{2.2.4} $\ad{T} = \sup \{n>0:
	 \|T_t\|_{1\to\infty} \leq Ct^{-n/2},\ t\geq1 \}$.
 \end{Prop}
 \begin{proof}
	 Set $d$ for the supremum.  Then for all $\eps>0$, there is 
	 $t_0>1$ s.t. $\|T_t\|_{1\to\infty}\leq t^{-(\ad{T}-\eps)/2}$, for 
	 all $t\geq t_0$, and, by previous proposition, $\ad{T}-\eps\leq 
	 d$.  Conversely $\|T_t\|_{1\to\infty}\leq t^{-(d-\eps)/2}$, for 
	 all $t\geq1$ implies $d-\eps\leq\ad{T}$.  
 \end{proof}\medskip
 
  \begin{rem} Varopoulos, Saloff-Coste and Coulhon call dimension at
 $\infty$ of the semigroup any of the numbers $n$ verifying the
 equivalent conditions of Theorem \ref{2.2.2}. Clearly such dimensions
 form a left half line, and the previous Proposition shows that
 $d_\infty(T)$ coincides with its upper bound.
 \end{rem}

 We want to give a formula for the computation of the asymptotic 
 dimension of a semigroup, in the special case of a semigroup of 
 integral operators with continuous kernel.  We need some preliminary 
 results. So let $X$ be a Hausdorff topological space, and $\mu$ a Borel 
 measure on it with $\supp\mu=X$.
 
 \begin{Lemma}\label{le:norm}
	 Let $K$ be an integral operator with kernel $k\in C(X\times X)\cap
	 L^{\infty}(X\times X,\cam\otimes\cam,\mu\otimes\mu; End(V))$, where 
	 $V$ is endowed with a scalar product.  
	 Then $\|K\|_{1\to\infty} = \sup_{x\in X} \|k(x,x)\|$.
 \end{Lemma}
 \begin{proof}
	We begin by proving that $\|K\|_{1\to\infty} = M := \sup_{x,y\in X} 
	\|k(x,y)\|$. Recall that
	\begin{equation*}
		\|K\|_{1\to\infty} = \sup_{f,g\in\O} \left\{ \left| 
		\int_{X}\int_{X} \langle f(x),k(x,y)g(y) \rangle 
		d\mu(x)d\mu(y) \right| \right\},
	\end{equation*}
	where $\O := \{f\in L^{1}(X,\cam,\mu;V), \|f\|=1 \}$.  Then it is 
	easy to see that $\|K\|_{1\to\infty}\leq M$.  For the reversed 
	inequality, let $\eps>0$, $(x_{\eps},y_{\eps})\in X\times X$ 
	be s.t. $M-\eps <\|k(x_{\eps},y_{\eps})\| \leq M$, and 
	$v_{\eps},w_{\eps}\in V$ be s.t. $\langle 
	v_{\eps},k(x_{\eps},y_{\eps})w_{\eps} \rangle \geq 
	\|k(x_{\eps},y_{\eps})\|-\eps$.  Let $A_{\eps},B_{\eps} \subset X$ 
	be open neighbourhoods of $x_{\eps}$, respectively $y_{\eps}$, of 
	finite measure s.t. $M-\eps <\|k(x,y)\| \leq M$, for any 
	$(x,y)\in A_{\eps}\times B_{\eps}$, and let $f_{\eps}(x) := 
	\frac{\chi_{A_{\eps}}(x)}{\mu(A_{\eps})}v_{\eps}$, $g_{\eps}(y) := 
	\frac{\chi_{B_{\eps}}(y)}{\mu(B_{\eps})}w_{\eps}$.  Then
	\begin{align*}
		& \left| \int_{X}\int_{X} \langle 
		f_{\eps}(x),k(x,y)g_{\eps}(y) \rangle d\mu(x)d\mu(y) \right| 
		\\
		& = \frac{1}{\mu(A_{\eps})\mu(B_{\eps})} 
		\left|\int_{A_{\eps}}d\mu(x)\int_{B_{\eps}}d\mu(y) \langle 
		v_{\eps},k(x,y)w_{\eps} \rangle \right| \\
		& \geq \frac{1}{\mu(A_{\eps})\mu(B_{\eps})} 
		\int_{A_{\eps}}d\mu(x)\int_{B_{\eps}}d\mu(y) \langle 
		v_{\eps},k(x_{\eps},y_{\eps})w_{\eps} \rangle + \\
		& - \frac{1}{\mu(A_{\eps})\mu(B_{\eps})} 
		\int_{A_{\eps}}d\mu(x)\int_{B_{\eps}}d\mu(y) \left|\langle 
		v_{\eps},[k(x,y)-k(x_{\eps},y_{\eps})]w_{\eps} \rangle \right| 
		\\
		& \geq \|k(x_{\eps},y_{\eps})\|-\eps - 
		\frac{1}{\mu(A_{\eps})\mu(B_{\eps})} 
		\int_{A_{\eps}}d\mu(x)\int_{B_{\eps}}d\mu(y) 
		\|k(x,y)-k(x_{\eps},y_{\eps})\| \\
		& \geq \|k(x_{\eps},y_{\eps})\| - 3\eps \geq M-4\eps.
	\end{align*}
	So that $\|K\|_{1\to\infty} = M$ follows.  Therefore to prove the 
	thesis it suffices to show that $M_{0}:= \sup_{x\in X} \|k(x,x)\| 
	= M$.  As $M_{0}\leq M$ is obvious, we show the opposite 
	inequality.  Let $\eps>0$, and $x_{\eps}, y_{\eps}\in X$, 
	$v_{\eps},w_{\eps} \in V$, $A_{\eps}, B_{\eps} \subset X$ be as 
	above.  Then
	\begin{align*}
		|\langle f_{\eps},Kf_{\eps} \rangle| & = 
		\frac{1}{\mu(A_{\eps})^{2}} 
		\left|\int_{A_{\eps}}d\mu(x)\int_{A_{\eps}}d\mu(y) \langle 
		v_{\eps},k(x,y)v_{\eps} \rangle \right| \\
		& \leq \frac{1}{\mu(A_{\eps})^{2}} 
		\int_{A_{\eps}}d\mu(x)\int_{A_{\eps}}d\mu(y) \langle 
		v_{\eps},k(x_{\eps},x_{\eps})v_{\eps} \rangle + \\
		& + \frac{1}{\mu(A_{\eps})^{2}} 
		\int_{A_{\eps}}d\mu(x)\int_{A_{\eps}}d\mu(y) \left|\langle 
		v_{\eps},[k(x,y)-k(x_{\eps},x_{\eps})]v_{\eps} \rangle \right| \\
		& \leq \|k(x_{\eps},x_{\eps})\| + \sup_{x,y \in 
		A_{\eps}}\|k(x,y)-k(x_{\eps},x_{\eps})\| \leq M_{0} + \eps,
	\end{align*}
	where the last inequality follows from the continuity of $k$, if 
	we choose $A_{\eps}$ small enough.  Analogously $|\langle 
	g_{\eps},Kg_{\eps} \rangle| \leq M_{0} + 2\eps$.  Then using the 
	estimates proved above and Cauchy-Schwarz inequality, we obtain
	\begin{align*}
		M- 4\eps & \leq |\langle f_{\eps},Kg_{\eps} \rangle| \\
		& \leq |\langle f_{\eps},Kf_{\eps} \rangle|^{1/2} |\langle 
		g_{\eps},Kg_{\eps} \rangle|^{1/2} \\
		& \leq M_{0}+ 2\eps,
	\end{align*}
	and from the arbitrariness of $\eps$ we get the thesis.
 \end{proof}
   
 \begin{Lemma}\label{le:pos}
 	If $\langle f,Kf \rangle \geq 0$ for any $f\in L^{2}(X,\cam,\mu;V)$, 
 	then $\langle v,k(x,x)v \rangle \geq 0$ for any $x\in X$, $v\in V$.
 \end{Lemma}
 \begin{proof}
	Assume on the contrary that there are $x_{0}\in X$, $v\in V$ s.t. 
	$\langle v,k(x_{0},x_{0})v \rangle < 0$.  Then, by continuity, 
	there is an open neighbourhood $U$ of $x_{0}$ s.t. $\re \langle 
	v,k(x,y)v \rangle$ $< 0$, for any $x,y \in U$. Let $f_{U}(x) := 
	\frac{\chi_{U}(x)}{\mu(U)}v$, so that
	\begin{align*}
		0\leq \langle f_{U},Kf_{U} \rangle & = \re \langle f_{U},Kf_{U} 
		\rangle\\
		& = \frac{1}{\mu(U)^{2}} 
		\int_{U}d\mu(x) \int_{U} d\mu(y) \re \langle v,k(x,y)v \rangle < 0
	\end{align*}
	which is absurd.
 \end{proof}
 
 \begin{Thm}\label{thm:asdimp_pos_op}
	 Let $X$ be a Hausdorff space, $\mu$ a Borel measure on it, 
	 $T_{t}: L^1(X,\cam,\m;V)\to L^\infty(X,\cam,\m;V)$ a semigroup of 
	 integral operators with continuous kernels $k(t,x,y)$, satisfying the 
	 hypotheses of Proposition \ref{2.2.3}, and assume $T_{t}$ is a 
	 positive bounded operator on $L^{2}(X,\cam,\m;V)$.  Then 
	 \begin{equation*}
		\ad{T} = \liminf_{t\to\infty} \frac{-2 \log(\sup_{x\in 
		X}Tr(k(t,x,x)))}{\log t}.
	 \end{equation*}
 \end{Thm}
 \begin{proof}
	 In the following we use the notation $f\pesci g$, where $f,g : 
	 [0,\infty)\to [0,\infty)$ to say that there are $t_{0}>0$, $C>0$ s.t. 
	 $C^{-1}\leq \frac{f(t)}{g(t)} \leq C$, for any $t\geq t_{0}$.  As 
	 $\|\cdot\|_{\infty}$ and $\|\cdot\|_{1}$ are equivalent on 
	 $End(V)$, and using Lemmas \ref{le:norm}, \ref{le:pos}, we get
 \begin{align*}
 	\|T_{t}\|_{1\to\infty} & = \sup_{x\in X} \|k(t,x,x)\| \\ 
	& \pesci \sup_{x\in X} Tr(|k(t,x,x)|) = \sup_{x\in X} Tr(k(t,x,x)). 
 \end{align*}
 Therefore 
 \begin{align*}
	\ad{T} & = \liminf_{t\to\infty} \frac{-2\log 
	\|T_{t}\|_{1\to\infty} }{ \log t} \\
	& = \liminf_{t\to\infty} \frac{-2\log (\sup_{x\in X} Tr(k(t,x,x))) 
	}{ \log t} .
 \end{align*}
 \end{proof}\medskip

  Using these results we can show the relation between the asymptotic 
  dimension of the heat kernel semigroup and the Novikov-Shubin 
  numbers.

 \begin{Cor}\label{asdim=NSinv}
	 Let $M$ be a complete connected Riemannian $n$-manifold, and $G$ 
	 an infinite discrete group of isometries of $M$, acting freely 
	 and properly discontinuously, and with $X:=M/ G$ a compact 
	 manifold.  Then $\ad{\e{-t\D_{k}}} = \underline\a_{k}(M,G)$.
 \end{Cor}
 \begin{proof}
	 In the following we use the notation $f\pesci g$, as in the proof 
	 of Theorem \ref{thm:asdimp_pos_op}.  Let us denote by 
	 $H_{k}(t,x,y)$ the kernel of the integral operator 
	 $\e{-t\D_{k}}$, and observe that
	 \begin{equation*}
		 \sup_{x\in M} Tr(H_{k}(t,x,x)) = \sup_{x\in\cf} 
		 Tr(H_{k}(t,x,x)) \pesci \inf_{x\in\cf} Tr(H_{k}(t,x,x)),
	 \end{equation*}
	 where the last relation follows from the fact that 
	 $\overline{\cf}$ is compact and $Tr(H_{k}(t,x,x))$ $>0$.  Therefore
	 $$
	 \sup_{x\in M} Tr(H_{k}(t,x,x)) \pesci \int_{\cf} 
	 Tr(H_{k}(t,x,x))dvol(x) = \th_{k}(t)-b_{k}.
	 $$
	 Then, using Theorem \ref{thm:asdimp_pos_op}, we get 
	 \begin{align*}
		 \ad{T} & = \liminf_{t\to\infty} \frac{-2\log 
		 \|\e{-t\D_{k}}\|_{1\to\infty} }{ \log t} \\
		 & = \liminf_{t\to\infty} \frac{-2\log (\th_{k}(t)-b_{k}) }{ 
		 \log t} = \underline\a_{k}(M,G).
	 \end{align*}
 \end{proof}\medskip

\subsection{Comparison between the algebras associated to a covering 
manifold and a general open manifold}
 
 In this subsection we study the relation between the von 
 Neumann algebra of $G$-invariant operators considered here, and the 
 C$^{*}$-algebra of almost local operators considered in \cite{GI2}, 
 namely the norm closure of the finite propagation operators, and the 
 traces on these algebras.

 \begin{Prop} 
	 Let us denote by $\Aloc_k$ the C$^*$-algebra of almost local 
	 operators on $k$-forms.  Then $\Aloc_k\cap\Ainv_k$ is weakly 
	 dense in $\Ainv_k$.
 \end{Prop}
 \begin{proof} 
	 First we choose a fundamental domain $\cf$ in $M$ and denote by 
	 $e_g$ the projection given by the multiplication operator by the 
	 characteristic function of $g\cf$, $g\in G$.  Then denote by $G_n$ 
	 the ball of radius $n$ in $G$, namely the set of elements which 
	 can be written as words of length $\leq n$ in terms 
	 of a prescribed set of generators for $G$. \\
	 For any selfadjoint operator $a$ acting on $L^2\La^k(M)$ and any 
	 $n\in\bn$ set
	 $$
	 a_n:=\sum_{g^{-1}h\in G_n}e_gae_h,
	 $$
	 and note that $a_n$ has finite propagation, hence it belongs to $\Aloc_k$. \\
	 Observe then that if $a$ is bounded, $a_n$ is bounded too. Indeed
	 \begin{align*}
		 (x,a_nx) &=\sum_{g^{-1}h\in G_n}(x,e_gae_h x) \leq \|a\| 
		 \sum_{g^{-1}h\in G_n}(e_gx,e_h x)\cr &\leq \|a\| 
		 \sum_{g^{-1}h\in G_n}\frac{\|e_gx\|^2+\|e_hx\|^2}2 = \|a\| 
		 \|x\|\; \#(G_n).
	 \end{align*}
	 Also, if $a$ is periodic, $a_n$ is periodic too. Indeed, for any 
	 $\g\in G$,
	 $$
	 L_\g a_n=\sum_{g^{-1}h\in G_n}L_\g e_gae_h=\sum_{g^{-1}h\in 
	 G_n}e_{\g g}ae_{\g h}L_\g =a_nL_\g.
	 $$
	 Since $a_n$ converges weakly to $a$, the thesis follows.
 \end{proof}
 
 Consider now the case that $G$ is amenable, and the corresponding (regular)
 exhaustion $\ck$ on $M$ (see \cite{Roe}). Denoting with $Tr_{\ck}$ 
 the trace on $\ca_{k}$ introduced in \cite{GI2}, the following holds
 
 \begin{Cor}
 	$Tr_{G}$ and $Tr_{\ck}$ coincide on $\Aloc_k\cap\Ainv_k$, hence 
 	$Tr_{G}$ is uniquely determined by $Tr_{\ck}$.
 \end{Cor}
 \begin{proof}
 	In this case $Tr_{\ck}$ is given, for $T\in \Aloc_k\cap\Ainv_k$, by 
	$$
	Tr_{\ck}(T) \equiv \frac{Tr(e_{1}Te_{1})}{vol(\cf)} = Tr_{G}(T),
	$$
	where $1\in G$ is the identity element, and we have chosen $vol(\cf)=1$. 
 \end{proof}


\end{document}